\documentclass[12pt]{article}
\usepackage{algorithm}
\usepackage{algorithmic}
\usepackage{mathrsfs}
\usepackage[justification=centering]{caption}
\usepackage{booktabs}
\usepackage{multirow}
\usepackage{array}
\usepackage[numbers,sort&compress]{natbib}
\usepackage[colorlinks,
            linkcolor=blue,
            anchorcolor=green,
            citecolor=magenta
            ]{hyperref}
\usepackage{mathrsfs,amssymb,amsfonts}
\usepackage{graphicx}
\usepackage{amsmath,amssymb,latexsym,color}
\usepackage[mathscr]{eucal}
\usepackage[numbers]{natbib}
\usepackage{verbatim}
\newtheorem{thm}{Theorem}[section]
\newtheorem{definition}[thm]{Definition}

\newtheorem{lemma}[thm]{Lemma}

\newtheorem{corollary}[thm]{Corollary}
\newtheorem{example}{Example}[section]
\newtheorem{remark}{Remark}[section]
\newtheorem{pro}{Problem}[section]
\newtheorem{observation}{Observation}[section]

\newcommand{\proof}{{\it Proof.\quad}}
\newcommand{\qed}{\hfill\Box\medskip}
\usepackage{CJK}
\usepackage{setspace}

\begin{document}

\title{\bf Biharmonic distance of graphs
}

\author{
Yulong Wei\textsuperscript{1}
\quad
Rong-hua Li\textsuperscript{2}\footnote{\footnotesize Corresponding authors.\newline
\indent\indent{\em E-mail address:} weiyulong@tyut.edu.cn (Y. Wei), lironghuabit@126.com (R-H. Li), yangweihua@tyut.edu.cn(W. Yang).}
\quad
Weihua Yang\textsuperscript{1}
\\
\\{\footnotesize
\textsuperscript{1}\em Department of Mathematics, Taiyuan University of Technology, Taiyuan, 030024, China}
\\{\footnotesize \textsuperscript{2}\em School of Computer Science {\rm \&} Technology, Beijing Institute of Technology, Beijing, 100081, China}}
\date{}
\maketitle

\setlength{\baselineskip}{24pt}

\noindent {\bf Abstract}\quad Lipman et al. [ACM Transactions on Graphics 29 (3) (2010), 1--11] introduced the concept of biharmonic distance to measure the distances between pairs of points on a 3D surface. Biharmonic distance has some advantages over resistance distance and geodesic distance in some realistic contexts. Nevertheless, limited work has been done on the biharmonic distance in the discrete case. In this paper, we give some characterizations of the biharmonic distance of a graph. Some basic mathematical properties of biharmonic distance and biharmonic index are established.


\vskip0.6cm

\section{Introduction}
Measuring the distances between pairs of vertices of a graph is a classical problem and thus several distance measures have
previously been proposed and are commonly used in many fields.
For example, the geodesic distance is using the length of the shortest path to measure distance which is intuitive and useful.
The geodesic distance of any two vertices of a graph can be computed by the
breadth-first search algorithm. The diameter of a graph is the maximal geodesic distance between vertices. Network topology is always represented by a graph, where vertices represent
processors and edges represent links between processors. The diameter is an important parameter for interconnection networks.
As another example, resistance distance,
which originated in electrical circuit theory, has played a prominent role in circuit theory \cite{Doy}, chemistry \cite{Kle,Peng},
combinatorial matrix theory \cite{Bap,Yang} and spectral graph theory \cite{CheH}. Moreover, resistance distance has extensive applications ranging from quantifying biological structures \cite{Peng},
distributed control systems \cite{Bar} and power grid systems \cite{Ste}.

The concept of biharmonic distance was first proposed in \cite{Lip} as a measure of distance between two points on a curve surface.
It has some advantages over geodesic distance and resistance distance as a metric that incorporates both
local and global graph structure in some realistic contexts, e.g., computer graphics, geometric processing, and shape analysis. Moreover, it can be used to measure the robustness of the second order noisy consensus problem without leaders \cite{Bam}. In 2018, Yi et al. \cite{Yi} defined the biharmonic distance index (i.e., biharmonic Kirchhoff index \cite{YiY}, or biharmonic index \cite{Xu}) of a graph to further
describe the behavior of second-order consensus dynamics and established a connection between the biharmonic distance of a
graph and its second-order network coherence. Since direct computation of biharmonic distances is computationally infeasible for huge networks with millions of vertices, Zhang et al. \cite{Zhang} developed a nearly linear-time algorithm to approximate all diagonal entries of pseudoinverse
of the square of graph Laplacian matrix by the Johnson-Lindenstrauss lemma and Laplacian
solvers. However, limited work has been done on the biharmonic distance of a graph.

In this paper, we investigate some properties on the biharmonic distance and biharmonic index of an undirected connected simple graph. First, we give some characterizations of the biharmonic distance of a graph and obtain some bounds for it. Second, we reveal a relationship between biharmonic index and Kirchhoff index and determine the unique graph having the minimum
biharmonic index among the connected graphs with given order. As applications, we study the biharmonic distance of some special graphs, such as the complement of a graph, the Cartesian product of
two graphs and the Cayley graph of a finite abelian group. In particular, the combinatorial expression of biharmonic distance between any two vertices of hypercubes are determined.

The rest of this paper is organized as follows. Section \ref{B} introduces some terminology and preliminaries. Some characterizations of the biharmonic distance of a graph are given in Section \ref{C}. The biharmonic index of a graph is investigated in Section \ref{D}. Some applications are given in Section \ref{E}. Section \ref{F} concludes the paper and proposes some open problems.

\section{Terminology and preliminaries}\label{B}
Let $G$ be a connected graph with $V(G)=\{v_1,\ldots,v_n\}$ and $E(G)=\{e_1,\ldots,e_m\}$ and let $G^c$ be the complement of a graph $G$. The {\em neighborhood} $N_G(v)$ of a vertex $v$ in $G$ is the set of vertices adjacent to $v$. The degree of a vertex $v$ in $G$ is $|N_G(v)|$. A graph $G$ is $k$-{\em regular} if $|N_G(v)|=k$ for any $v\in V(G)$. We refer readers to \cite{Bon} for terminology and notation unless stated otherwise.

The adjacency matrix of the graph is denoted by $A$ and $D$ represents the diagonal matrix of vertex degrees.
The Laplacian matrix of $G$ is the matrix $D-A$, denoted by $L$ (or $L(G)$).
Throughout, $O$ denotes the all-$0$ matrix, $\mathbf{0}$ the all-$0$ vector and $\mathbf{1}$ the all-$1$ vector. Suppose that $\delta_u$ is the elementary unit vector with a $1$ in coordinate $u$. In the whole paper, we write $\delta_{uv}$ instead of $\delta_u-\delta_v$.

The following result shows that the spectrum of $L$ determines the number of connected components in a graph.
\begin{lemma}[\cite{Cve}]\label{M0}
The multiplicity of $0$ as an eigenvalue of $L$ is equal to the
number of connected components in $G$.
\end{lemma}
Thus, if $G$ is a connected graph, then the multiplicity of $0$ as an eigenvalue of $L$ is $1$.
In fact, the matrix $L$ is symmetric and positive semidefinite, which means $L$
has a spectral decomposition as
\begin{equation*}
L=\sum_{k=1}^n\lambda_kz_kz_k^{\top},
\end{equation*}
where $0=\lambda_1<\lambda_2\leq\cdots\leq\lambda_n$ are its $n$ positive eigenvalues,
and $z_1,\ldots,z_n$ are the corresponding mutually orthogonal unit eigenvectors.

The following two lemmas are important to obtain our main results.

\begin{lemma}[\cite{Das0}]\label{Da}
Let $G$ be a graph with $n$ vertices and $E(G)\neq\emptyset$. Then
$\lambda_2(G)=\cdots=\lambda_n(G)$ if and only if $G$ is isomorphic to $K_n$.
\end{lemma}

\begin{lemma}[Matrix-Tree Theorem \cite{Bap}]\label{MT}
Let $G$ be a graph. Then the cofactor of
any element of $L$ equals the number of spanning trees of $G$.
\end{lemma}

Let $B$ be an $m\times n$ matrix. A matrix $H$ of order $n\times m$ is said to be a \emph{$g$-inverse matrix} of $B$ if $BHB=B$. The \emph{Moore Penrose inverse} of an $n\times m$ real matrix $M$, denoted by $M^{+}$,
is the $m\times n$ real matrix that satisfies the following equations:
$$MM^{+}M=M, M^{+}MM^{+}=M^{+},(MM^{+})^{\top}=MM^{+}, (M^{+}M)^{\top}=M^{+}M. $$
We call the Moore Penrose inverse of $L$ {\em pseudoinverse Laplacian} and write $L^{2+}$
instead of $(L^{+})^2$. Similarly, the matrix $L^{+}$ also has a spectral decomposition as
\begin{equation}\label{1}
L^{+}=\sum_{k=2}^n\lambda_k^{-1}z_kz_k^{\top},
\end{equation}
where $\lambda_k$ is the $k$th smallest eigenvalue of the Laplacian matrix $L$ of a connected graph $G$ for $k=2,\ldots,n$ and $z_2,\ldots,z_n$
are the corresponding mutually orthogonal unit eigenvectors. Actually, the eigenvalues of $L^{+}$ are $0,\lambda_n^{-1},\lambda_{n-1}^{-1},\ldots,\lambda_2^{-1}$ \cite{Bap}.

The notion of biharmonic distance was first proposed in \cite{Lip} as a measure of distance between two points $u$ and $v$
on a curve surface:
$$d_B^2(u,v)=g_d(u,u)+g_d(v,v)-g_d(u,v),$$
where $g_d$ is the discrete Green's function of the discretized, normalized bilaplacian $\widetilde{L}^2$, equivalent to the pseudoinverse of $\widetilde{L}^2$, and $\widetilde{L}$ is the normalization of Laplacian $L$.
Fitch and Leonard \cite{Fit} defined the biharmonic distance between two vertices $i$ and $j$ in the graph
$G$, which they denoted $\gamma_{i,j}$, analogously without normalizing $L$:
$$\gamma_{i,j}=L_{ii}^{2+}+L_{jj}^{2+}-2L_{ij}^{2+}.$$
In this paper we adopt the definition in \cite{Yi}, which is a fusion of the above two conceptions, defined as follows.
\begin{definition}[\cite{Lip,Fit,Yi}]\label{D1}
For a graph $G$, the biharmonic distance $d_B(u,v)$ between two vertices $u$ and $v$ is defined by
\begin{equation}\label{2}
d_B^2(u,v)=L_{uu}^{2+}+L_{vv}^{2+}-2L_{uv}^{2+},
\end{equation}
where $L^{2+}_{ij}$ means the $(i,j)$-entry of $L^{2+}$.
\end{definition}

Biharmonic distance is a metric, which satisfies the following properties: non-negativity, nullity, symmetry, and triangle inequality \cite{Lip}.
The biharmonic index \cite{Yi}, as a distance-based graph invariant, is used for describing the behavior of second-order consensus dynamics. The concept of biharmonic index of graphs is defined as follows.
\begin{definition}[\cite{Yi,Xu}]\label{D3}
For a connected graph $G$, the biharmonic index of $G$ is equal to
$$\frac{1}{2}\sum_{u\in V(G)}\sum_{v\in V(G)}d_B^2(u,v), $$
denoted by $B(G)$.
\end{definition}

\section{The biharmonic distance of a graph}\label{C}
In this section, we will give some characterizations of the biharmonic distance of a graph. Given a connected graph $G$, suppose that $\lambda_k$ is the $k$th smallest eigenvalue of the Laplacian matrix of $G$ for $k=1,2,\ldots,n$ and $z_1,z_2,\ldots,z_n$ are the corresponding mutually orthogonal real unit eigenvectors. Since $L\mathbf{1}=\mathbf{0}$, set $z_1=\frac{1}{\sqrt{n}}\mathbf{1}$. Let $B$ be an $m\times n$ matrix. If $S\subset\{1,\ldots,m\}$, $T\subset\{1,\ldots,n\}$, then $B[S|T]$ will
denote the submatrix of $B$ determined by the rows corresponding to $S$ and the columns
corresponding to $T$. The submatrix obtained by deleting the rows indexed in $S$ and the columns
indexed in $T$ will be denoted by $B(S|T)$. When $S=\{i\}$, $T=\{j\}$ are singletons, then $B(S|T)$ is denoted by $B(i|j)$.
The \emph{rank} of a matrix $B$ is denoted by $\textrm{r}(B)$.

\subsection{The spectral characterization}
\begin{thm}\label{m1}
Let $G$ be a connected graph with $n$ vertices. Then the biharmonic distance between two vertices $u$ and $v$ is
\begin{equation}\label{3}
d_B(u,v)=(\sum_{k=2}^{n}\lambda_k^{-2}(z_k(u)-z_k(v))^2)^{\frac{1}{2}}.
\end{equation}
\end{thm}
\proof  Note that $z_1,z_2,\ldots,z_n$ is a normal orthogonal basic of $\mathbf{R}^n$. Then by formulas (\ref{1}) and (\ref{2}),
\begin{eqnarray*}
d_B^2(u,v)&=&\delta_{uv}^{\top}L^{2+}\delta_{uv}\\
&=&\delta_{uv}^{\top}(L^{+}L^{+})\delta_{uv}\\
&=&\delta_{uv}^{\top}(\sum_{k=2}^n\lambda_k^{-1}z_kz_k^{\top})(\sum_{k=2}^n\lambda_k^{-1}z_kz_k^{\top})\delta_{uv}\\
&=&\delta_{uv}^{\top}(\sum_{k=2}^n\lambda_k^{-2}z_kz_k^{\top})\delta_{uv}\\
&=&\sum_{k=2}^n\lambda_k^{-2}(\delta_{uv}^{\top}z_kz_k^{\top}\delta_{uv})\\
&=&\sum_{k=2}^n\lambda_k^{-2}(z_k(u)-z_k(v))^2.
\end{eqnarray*}
Thus, we obtain the desired result. $\qed$

Let $S=\{k\mid 2\leq k\leq n\}$, $\sigma_2=\{k\mid 3\leq k\leq n ~\textmd{and}~ \lambda_k>\lambda_2\}$ and $\sigma_n=\{k\mid 2\leq k\leq n-1~\textmd{and}~\lambda_k<\lambda_n\}$. The eigenspace of $L$ corresponding to $\lambda_k$ is denoted by $V_{\lambda_k}$ and the orthogonal complement of $V_{\lambda_k}$ is denoted by $V_{\lambda_k}^{\bot}$ for $1\leq k\leq n$.
\begin{thm}\label{m2}
Let $G$ be a connected graph with $n$ vertices. Then the biharmonic distance between two distinct vertices $u$ and $v$ satisfies that
$$
\sqrt{2}\lambda_n^{-1}\leq d_B(u,v)\leq\sqrt{2}\lambda_2^{-1},
$$
and $d_B(u,v)=\sqrt{2}\lambda_j^{-1}$ if and only if $G\cong K_n$ or $\delta_{uv}\in V_{\lambda_k}^{\bot}$ for any $k\in\sigma_j$, where $j\in\{2,n\}$.
\end{thm}
\proof Since
\begin{eqnarray*}
&&\sum_{k=2}^{n}(z_k(u)-z_k(v))^2\\
&=&\sum_{k=2}^{n}z_k^2(u)+\sum_{k=2}^{n}z_k^2(v)-2\sum_{k=2}^{n}z_k(u)z_k(v)\\
&=&(\sum_{k=1}^{n}z_k^2(u)-z_1^2(u))+(\sum_{k=1}^{n}z_k^2(v)-z_1^2(v))-2(\sum_{k=1}^{n}z_k(u)z_k(v)-z_1(u)z_1(v))\\
&=&(1-\frac{1}{n})+(1-\frac{1}{n})-2(0-\frac{1}{\sqrt{n}}\times\frac{1}{\sqrt{n}})\\
&=&2,
\end{eqnarray*}
by formula (\ref{3}), we have
$$d_B(u,v)\leq\lambda_2^{-1}(\sum_{k=2}^{n}(z_k(u)-z_k(v))^2)^{\frac{1}{2}}=\sqrt{2}\lambda_2^{-1}.$$

Note that the Laplacian eigenvalues of $K_n$ are $n$ with multiplicities $n-1$ and $0$ with multiplicity $1$. Then $d_B(u,v)=(\sum_{k=2}^{n}n^{-2}(z_k(u)-z_k(v))^2)^{\frac{1}{2}}=\dfrac{\sqrt{2}}{n}$ for any $u$, $v\in V(K_n)$.
Therefore, if $G\cong K_n$, then the equality holds.

If $\delta_{uv}\in V_{\lambda_k}^{\bot}$ for any $k\in\sigma_2$, then $z_k(u)-z_k(v)=\delta_{uv}^{\top}z_k=0$ for any $z_k\in V_{\lambda_k}$ and $k\in\sigma_2$. Thus,
\begin{eqnarray*}
d_B^2(u,v)&=&\sum_{k=2}^{n}\lambda_k^{-2}(z_k(u)-z_k(v))^2\\
&=&\sum_{k\in\sigma_2}\lambda_k^{-2}(z_k(u)-z_k(v))^2+\sum_{k\in S\setminus\sigma_2}\lambda_k^{-2}(z_k(u)-z_k(v))^2\\
&=&\sum_{k\in\sigma_2}(\lambda_k^{-2}\times 0)+\sum_{k\in S\setminus\sigma_2}\lambda_2^{-2}(z_k(u)-z_k(v))^2\\
&=&\lambda_2^{-2}\sum_{k=2}^{n}(z_k(u)-z_k(v))^2\\
&=&2\lambda_2^{-2}.
\end{eqnarray*}
Thus, the equality holds.

On the other hand, if $d_B(u,v)=\sqrt{2}\lambda_2^{-1}$ and $G$ is not isomorphic to $K_n$, then by Lemma \ref{Da}, $\sigma_2\neq\emptyset$.
Assume to the contrary that there exists an integer $k_1\in \sigma_2$ such that $\delta_{uv}\notin V_{\lambda_{k_1}}^{\bot}$. Then $z_{k_1}(u)-z_{k_1}(v)=\delta_{uv}^{\top}z_{k_1}\neq0$.
Thus,
\begin{eqnarray*}
d_B^2(u,v)&=&\sum_{k=2}^{n}\lambda_k^{-2}(z_k(u)-z_k(v))^2\\
&=&\lambda_{k_1}^{-2}(z_{k_1}(u)-z_{k_1}(v))^2+\sum_{k\in S\setminus\{k_1\}}\lambda_k^{-2}(z_k(u)-z_k(v))^2\\
&<&\lambda_2^{-2}(z_{k_1}(u)-z_{k_1}(v))^2+\sum_{k\in S\setminus\{k_1\}}\lambda_2^{-2}(z_k(u)-z_k(v))^2\\
&=&\lambda_2^{-2}\sum_{k=2}^{n}(z_k(u)-z_k(v))^2\\
&=&2\lambda_2^{-2},
\end{eqnarray*}
a contradiction.

By the similar arguments, we have
$d_B(u,v)\geq\sqrt{2}\lambda_n^{-1}$ and equality holds if and only if $G\cong K_n$ or $\delta_{uv}\in V_{\lambda_k}^{\bot}$ for any $k\in\sigma_n$.

As mentioned above, we complete this proof. $\qed$

\begin{figure}[hptb]
  \centering
  \includegraphics[width=8cm]{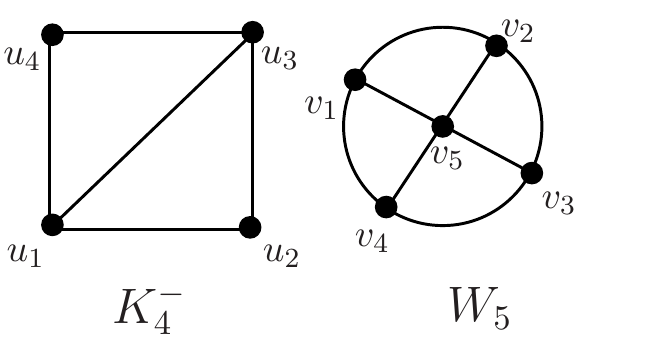}\\
  \caption{The structures of $K_4^-$ and $W_5$.
}\label{f1}
\end{figure}

\begin{example}\label{SEC}
Let $K_4^-$ be the graph obtained from $K_4$ by deleting an edge and let $W_5$ be a wheel graph with $5$ vertices (see Figure \ref{f1}). The eigenvalues of the Laplacian matrix of $K_4^-$ are $0$, $2$, $4$, $4$. Then $\sigma_4=\{2\}$.
Suppose $x=(x_{u_1},\ldots,x_{u_4})^{\top}\in \mathbf{R}^4$. Solve the system of linear equations $(2I-L(K_4^-))x=\mathbf{0}$ and obtain that $x_{u_1}=x_{u_3}=0$. Then $d_B(u_1,u_3)=\dfrac{\sqrt{2}}{4}$ in $K_4^-$. The eigenvalues of the Laplacian matrix of $W_5$ are $0$, $3$, $3$, $5$, $5$. Then $\sigma_2=\{4,5\}$.
Suppose $y=(y_{v_1},\ldots,y_{v_5})^{\top}\in \mathbf{R}^5$. Solve the system of linear equations $(5I-L(W_5))x=\mathbf{0}$ and obtain that $y_{v_2}=y_{v_4}$. Then $d_B(v_2,v_4)=\dfrac{\sqrt{2}}{3}$ in $W_5$.
\end{example}

\subsection{The determinant characterization}
Let $L$ be the Laplacian matrix of a connected graph $G$ with $n$ vertices. In this section, we first give a formula for biharmonic distance of $G$ in terms of the $g$-inverse matrix of $L^2$ and then by this formula, an expression for biharmonic distance between any two distinct vertices of $G$ is established in terms of an $(n-2)\times (n-2)$ principal minor of $L^2$ and the number of spanning trees of $G$.
\begin{lemma}\label{L4}
Let $G$ be a connected graph with $n$ vertices. If $u, v\in V(G)$, then
$$d_B^2(u,v)=\delta_{uv}^{\top}H\delta_{uv},$$
where $L$ is the Laplacian matrix of $G$ and $H$ is a $g$-inverse matrix of $L^2$.
\end{lemma}
\proof
\noindent $\mathbf{Claim~1.}$~~$L^{2+}=(L^2)^+$.

By formula (\ref{1}), $L^{2+}=L^{+}L^{+}=(\sum_{k=2}^n\lambda_k^{-1}z_kz_k^{\top})^2=\sum_{i=2}^n\lambda_k^{-2}z_kz_k^{\top}=(L^2)^+$.

\noindent $\mathbf{Claim~2.}$~~$\delta_{uv}^{\top}H_1\delta_{uv}=\delta_{uv}^{\top}H_2\delta_{uv}$ for any two $g$-inverse matrices $H_1$ and $H_2$ of $L^2$.

We use Col$(B)$ to denote the column space of a matrix $B$. Since $L$ is a symmetric matrix, we have $\textrm{r}(L^2)=\textrm{r}(L^{\top}L)=\textrm{r}(L)$. Thus, $\dim$(Col($L^2))=\dim$(Col$(L))$. Note that Col$(L^2)$$\subseteq$Col$(L)$. Then Col$(L^2)$$=$Col$(L)$. Since $G$ is a connected graph, by Lemma \ref{M0}, $\textrm{r}(L)=n-1$. Since $L\mathbf{1}=\mathbf{0}$ and $\delta_{uv}^{\top}\mathbf{1}=0$, $\delta_{uv}\in$Col$(L)$. Then $\delta_{uv}$$\in$Col$(L^2)$. Therefore, there exists a vector $z$ such that $L^2z=\delta_{uv}$. Thus, considering $H_1$, $H_2$ any two $g$-inverse matrices of $L^2$,
\begin{eqnarray*}
\delta_{uv}^{\top}H_1\delta_{uv}-\delta_{uv}^{\top}H_2\delta_{uv}&=&\delta_{uv}^{\top}(H_1-H_2)\delta_{uv}\\
&=&(L^2z)^{\top}(H_1-H_2)(L^2z)\\
&=&z^{\top}((L^2)^{\top}(H_1-H_2)L^2)z\\
&=&z^{\top}(L^2H_1L^2-L^2H_2L^2)z\\
&=&z^{\top}(L^2-L^2)z\\
&=&0.
\end{eqnarray*}

By Claims 1 and 2, we have $d_B^2(u,v)=\delta_{uv}^{\top}L^{2+}\delta_{uv}=\delta_{uv}^{\top}(L^2)^+\delta_{uv}=\delta_{uv}^{\top}H\delta_{uv}$, where $H$ is a $g$-inverse matrix of $L^2$. $\qed$

\begin{thm}\label{m3}
Let $G$ be a connected graph with $n$ vertices. Then the biharmonic distance between two distinct vertices $u$ and $v$ is
\begin{equation*}
d_B(u,v)=\dfrac{\sqrt{\det(L^2(u,v|u,v))}}{\sqrt{n}\tau(G)},
\end{equation*}
where $L$ is the Laplacian matrix of $G$ and $\tau(G)$ is the number of spanning trees of $G$.
\end{thm}
\proof Let $C=L^2(v|v)$. Then by Cauchy-Binet Formula and Lemma \ref{MT}, we have
\begin{eqnarray*}
\det C&=&\det(L^2(v|v))\\
&=&\sum_{S\subseteq V(G),|S|=n-1}(\det(L[V(G)\setminus\{v\}|S]))^2\\
&=&n\tau^2(G)\\
&>&0.
\end{eqnarray*}
We construct the following matrix $H$: In $L^2$, replace entries in the $v$th row and column by zeros and replace $L^2(v|v)$ by
$C^{-1}$.

\noindent $\mathbf{Claim.}$~~$L^2HL^2=L^2$.

By the construction of $H$, there exists a permutation matrix $P$ such that
$P^{\top}HP=\begin{bmatrix} 0 & O \\ O & Q^{\top}C^{-1}Q \end{bmatrix}$ and
$P^{\top}L^2P=\begin{bmatrix} b_{11} & b^{\top} \\ b & Q^{\top}CQ \end{bmatrix}$,
where $Q$ is a real $(n-1)\times(n-1)$ matrix such that $Q^{\top}Q=I$. Then
\begin{eqnarray*}
P^{\top}L^2HL^2P&=&(P^{\top}L^2P)(P^{\top}HP)(P^{\top}L^2P)\\
&=&\begin{bmatrix} b_{11} & b^{\top} \\ b & Q^{\top}CQ \end{bmatrix}\begin{bmatrix} 0 & O \\ O & Q^{\top}C^{-1}Q \end{bmatrix}\begin{bmatrix} b_{11} & b^{\top} \\ b & Q^{\top}CQ \end{bmatrix}\\
&=&\begin{bmatrix} 0 & b^{\top}Q^{\top}C^{-1}Q \\ O & I \end{bmatrix}\begin{bmatrix} b_{11} & b^{\top} \\ b & Q^{\top}CQ \end{bmatrix}\\
&=&\begin{bmatrix} b^{\top}Q^{\top}C^{-1}Qb & b^{\top} \\ b & Q^{\top}CQ \end{bmatrix}.
\end{eqnarray*}
Note that $\textrm{r}(L^2)=\textrm{r}(L)=n-1=\textrm{r}(C)$.
Then
\begin{eqnarray*}
\textrm{r}(Q^{\top}CQ)&=&\textrm{r}(P^{\top}L^2P)\\
&=&\textrm{r}\left(\begin{bmatrix} b_{11} & b^{\top} \\ b & Q^{\top}CQ \end{bmatrix}\right)\\
&=&\textrm{r}\left(\begin{bmatrix} b_{11}-b^{\top}(Q^{\top}CQ)^{-1}b & O \\ b & Q^{\top}CQ \end{bmatrix}\right)\\
&=&\textrm{r}\left(\begin{bmatrix} b_{11}-b^{\top}(Q^{\top}CQ)^{-1}b & O \\ O & Q^{\top}CQ \end{bmatrix}\right)\\
&=&\textrm{r}(Q^{\top}CQ)+\textrm{r}([b_{11}-b^{\top}(Q^{\top}CQ)^{-1}b]).
\end{eqnarray*}
Therefore, $\textrm{r}([b_{11}-b^{\top}(Q^{\top}CQ)^{-1}b])=0$ and $b_{11}=b^{\top}(Q^{\top}CQ)^{-1}b$. Thus, $P^{\top}L^2HL^2P=P^{\top}L^2P$ and so $L^2HL^2=L^2$.

According to the above claim, $H$ is a $g$-inverse of $L^2$.
By Lemma \ref{L4},
\begin{eqnarray*}
d_B^2(u,v)&=&\delta_{uv}^{\top}H\delta_{uv}\\
&=&h_{uu}+h_{vv}-h_{uv}-h_{vu}\\
&=&h_{uu}\\
&=&(C^{-1})_{uu}\\
&=&\dfrac{\det(L^2(u,v|u,v))}{\det(L^2(v|v))}\\
&=&\dfrac{\det(L^2(u,v|u,v))}{n\tau^2(G)}.
\end{eqnarray*}
This completes the proof of Theorem \ref{m3}. $\qed$

\subsection{The system of linear equations characterization}
In this section, we will characterize biharmonic distance between any two distinct vertices of a connected graph from the perspective of the system of linear equations.
Given a real vector $f\in\mathbf{R}^n$, let $\|f\|_2=\sqrt{f^{\top}f}$. We call $\|f\|_2$ the \emph{norm} of $f$.
\begin{thm}\label{m4}
Let $G$ be a connected graph with $n$ vertices. If $f$ is a solution of the equation $Lx=\delta_{uv}$ with minimum norm, then the biharmonic distance between two vertices $u$ and $v$ is
\begin{equation*}\label{7}
d_B(u,v)=\|f\|_2.
\end{equation*}
\end{thm}
\proof If $f$ is a solution of the equation $Lx=\delta_{uv}$ with minimum norm, then $f=L^+\delta_{uv}$. Thus,
\begin{eqnarray*}
d_B^2(u,v)&=&\delta_{uv}^{\top}L^{2+}\delta_{uv}\\
&=&(L^{+}\delta_{uv})^{\top}(L^{+}\delta_{uv})\\
&=&f^{\top}f\\
&=&\|f\|_2^2.
\end{eqnarray*}
Therefore, $d_B(u,v)=\|f\|_2$. $\qed$

\section{The biharmonic index of a graph}\label{D}
In this section, we investigate the biharmonic index $B(G)$ of a connected graph $G$. We first derive an expression for the biharmonic index of a graph $G$ in terms of the Laplacian eigenvalues of $G$. Then we reveal a relationship between biharmonic index and Kirchhoff index. Moreover, we study the variation of $B(G)$ when modifications on the graph are made, in particular when an edge is added. A lower bound for $B(G)$ is given and the equality case is studied.
\subsection{The spectral characterization}
\begin{thm}\label{m5}
Let $G$ be a connected graph with $n$ vertices and $\lambda_i$ the $i$th smallest eigenvalue of the Laplacian matrix of $G$. Then
\begin{equation*}\label{8}
B(G)=n\sum_{i=2}^{n}\dfrac{1}{\lambda_i^2}.
\end{equation*}
\end{thm}
\proof  Since
$$
L^{2+}\mathbf{1}=L^{+}L^{+}\mathbf{1}=L^{+}L^{+}LL^{+}\mathbf{1}=L^{+}L^{+}L^{+}L\mathbf{1}=\mathbf{0},
$$
we have $\sum_{v\in V(G)}L_{uv}^{2+}=0$ for any $u\in V(G)$. Then
\begin{eqnarray*}
B(G)&=&\dfrac{1}{2}\sum_{u\in V(G)}\sum_{v\in V(G)}d_B^2(u,v)\\
&=&\dfrac{1}{2}\sum_{u\in V(G)}\sum_{v\in V(G)}(L_{uu}^{2+}+L_{vv}^{2+}-2L_{uv}^{2+})\\
&=&n\sum_{u\in V(G)}L_{uu}^{2+}\\
&=&n~\text{trace}(L^{2+})\\
&=&n\sum_{i=2}^{n}\dfrac{1}{\lambda_i^2}.
\end{eqnarray*}
This completes the proof of Theorem \ref{m5}. $\qed$

\subsection{A relationship to Kirchhoff index}
For a graph $G$, its Kirchhoff index $Kf(G)$ \cite{Kle} is defined as the sum of resistance distance over all its vertex pairs.
That is, $Kf(G)=\sum_{i<j}r_{ij}$, where $r_{ij}$ is the resistance distance between the vertices $v_i$ and $v_j$ in $G$. This index can be
expressed via the non-zero Laplacian eigenvalues as follows.
\begin{lemma}[\cite{Gut}]\label{GM}
If $G$ be a connected graph with $n$ vertices, then
$$Kf(G)=n\sum_{i=2}^n\frac{1}{\lambda_i}.$$
\end{lemma}

Now, we establish a relationship between biharmonic index and Kirchhoff index in the following discussion.
\begin{thm}\label{BRK}
Let $G$ be a connected graph with $n$ vertices. Then for $n\geq2$,
$$B(G)\geq \dfrac{Kf(G)^2}{n(n-1)},$$
where the equality holds if and only if $G$ is isomorphic to $K_n$.
\end{thm}
\proof
By Theorem \ref{m5}, Cauchy-Schwarz inequality and Lemma \ref{GM}, we have
\begin{eqnarray*}
B(G)&=&n\sum_{i=2}^{n}\dfrac{1}{\lambda_i^2}\\
&=&\frac{n}{n-1}\cdot(n-1)\sum_{i=2}^{n}\dfrac{1}{\lambda_i^2}\\
&\geq&\frac{n}{n-1}(\sum_{i=2}^{n}\dfrac{1}{\lambda_i})^2\\
&=&\frac{n}{n-1}(\dfrac{Kf(G)}{n})^2\\
&=&\dfrac{Kf(G)^2}{n(n-1)}.
\end{eqnarray*}
Note that the equality holds if and only if $\lambda_2=\cdots=\lambda_n$. By Lemma \ref{Da}, we obtain the desired results. $\qed$

\subsection{The biharmonic index of connected graphs with given order}
\begin{lemma}[\cite{Cve}]\label{pur}
If $e$ is an edge of the graph $G$ and $G'=G-e$, then
$$0=\lambda_1(G')=\lambda_1(G)\leq\lambda_2(G')\leq\lambda_2(G)\leq\cdots\leq\lambda_{n-1}(G)\leq\lambda_n(G')\leq\lambda_n(G),$$
where $\lambda_1(G),\ldots,\lambda_n(G)$ and $\lambda_1(G'),\ldots,\lambda_n(G')$ are the Laplacian eigenvalues of $G$ and $G'$, respectively.
\end{lemma}

\begin{thm}\label{m6}
Let $G$ be a connected graph with $n$ vertices and $e\in E(G^c)$. Then
$$B(G+e)<B(G).$$
\end{thm}
\proof Note that $\sum_{i=2}^{n}\lambda_i(G)=\textrm{trace}~L(G)=2|E(G)|$ and $\sum_{i=2}^{n}\lambda_i(G+e)=2(|E(G)|+1)$. Then $\sum_{i=2}^{n}\lambda_i(G+e)-(\sum_{i=2}^{n}\lambda_i(G))=2$. By Lemma \ref{pur}, there exists some $j$ such that $\lambda_j(G+e)>\lambda_j(G)$. By Theorem \ref{m5}, we have
\begin{eqnarray*}
B(G+e)&=&n\sum_{i=2}^{n}\dfrac{1}{\lambda_i^2(G+e)}\\
&<&n\sum_{i=2}^{n}\dfrac{1}{\lambda_i^2(G)}\\
&=&B(G).
\end{eqnarray*}
This completes the proof of Theorem \ref{m6}. $\qed$

\begin{thm}\label{m7}
Let $G$ be a connected graph with $n$ vertices. Then
$$B(G)\geq \dfrac{n-1}{n},$$
where the equality holds if and only if $G$ is isomorphic to $K_n$.
\end{thm}
\proof
By the proof of Theorem \ref{m2}, $d_B(u,v)=\dfrac{\sqrt{2}}{n}$ for any two distinct vertices $u$, $v\in V(K_n)$. Therefore, by Theorem \ref{m6}, we have $$B(G)\geq B(K_n)=(\dfrac{\sqrt{2}}{n})^2\times \dbinom{n}{2}=\dfrac{n-1}{n}.$$
Suppose $B(G)=\dfrac{n-1}{n}$. Assume to the contrary that $G$ is not a complete graph and $e\in E(G^c)$. Then by Theorem \ref{m6}, we have $\dfrac{n-1}{n}\leq B(G+e)<B(G)=\dfrac{n-1}{n}$, a contradiction.
Thus, $B(G)=\dfrac{n-1}{n}$ if and only if $G$ is isomorphic to $K_n$.
$\qed$

\begin{remark}
Various bounds on the sum of powers of the Laplacian eigenvalues of graphs have been established
in terms of the number of vertices, number of edges, maximum degree, clique number, independence number, matching number or the number of spanning trees of graphs. We refer readers to \cite{CheD,Das,Liu,Zhou}. By Theorem \ref{m5}, other lower bounds on the biharmonic index of connected graphs can be obtained via the results in \cite{CheD,Das,Liu,Zhou}.
\end{remark}

\section{Applications to some special graphs}\label{E}
In this section, as applications of Theorem \ref{m1}, we investigate
the biharmonic distance of some special graphs including the complement of a graph,
the Cartesian product of two graphs and the Cayley graph of a finite abelian group. In
particular, the combinatorial expression of biharmonic distance between any two vertices
of hypercubes are determined.
\subsection{The complement of a graph}
Suppose that $\lambda_{1}\leq\cdots\leq\lambda_{n}$ are the eigenvalues of the Laplacian matrix of $G$ and $z_1,\ldots,z_n$ are the corresponding mutually orthogonal real unit eigenvectors.
\begin{thm}\label{MM}
Let $G^c$ be the complement of a graph $G$ with $n$ vertices. If $G^c$ is connected, then the biharmonic distance between two vertices $u$ and $v$ in $G^c$ is
$$d_B(u,v)=(\sum_{k=2}^{n}(n-\lambda_k)^{-2}(z_k(u)-z_k(v))^2)^{\frac{1}{2}}.$$
\end{thm}
\proof Note that $L(G^c)=nI-L(G)-J$, where $J$ denotes the all-$1$ matrix. Then $L(G^c)z_1=\mathbf{0}$ and
$L(G^c)z_k=(n-\lambda_{k})z_k$ for $k=2,\ldots,n$. Note that $\lambda_{n}+\lambda_{2}(G^c)=n$. Since $G^c$ is connected, $\lambda_{2}(G^c)>0$. Thus, $n-\lambda_k\geq n-\lambda_n=\lambda_{2}(G^c)>0$ for $2\leq k\leq n$. By Theorem \ref{m1}, we have $d_B(u,v)=(\sum_{k=2}^{n}(n-\lambda_k)^{-2}(z_k(u)-z_k(v))^2)^{\frac{1}{2}}$. $\qed$

\subsection{The Cartesian product of two graphs}
The \emph{Cartesian product} of two graphs $G_1$ and $G_2$ is the graph $G_1\square G_2$ whose vertex
set is the set $\{u_1u_2 | u_1\in V(G_1), u_2\in V(G_2)\}$, and two vertices $u_1u_2, v_1v_2$ are adjacent if $u_1=v_1$ and $\{u_2,v_2\}\in E(G_2)$ or if $u_2=v_2$ and $\{u_1,v_1\}\in E(G_1)$.
The \emph{Kronecker product} $A\otimes B$ of matrices $A=(a_{ij})_{m\times n}$
and $B=(b_{ij})_{p\times q}$ is the $mp\times nq$ matrix obtained from $A$ by replacing each
element $a_{ij}$ with the block $a_{ij}B$.

Suppose that $\lambda_{1}(G_i)\leq\cdots\leq\lambda_{n_i}(G_i)$ are the eigenvalues of the Laplacian matrix of $G_i$ and $z_1(G_i),\ldots,z_n(G_i)$ are the corresponding mutually orthogonal real unit eigenvectors for $i=1,2$.
\begin{lemma}[\cite{Fie}]\label{CP}
Let $G_1\square G_2$ be the Cartesian product of two graphs $G_1$ and $G_2$. Then the eigenvalues of the Laplacian matrix of $G_1\square G_2$ are
$\lambda_{i}(G_1)+\lambda_{j}(G_2)$ for $i=1,\ldots,n_1$ and $j=1,\ldots,n_2$ and the corresponding eigenvectors are $z_i(G_1)\otimes z_j(G_2)$ for $i=1,\ldots,n_1$ and $j=1,\ldots,n_2$.
\end{lemma}

By Theorem \ref{m1} and Lemma \ref{CP}, we immediately obtain the following result.
\begin{thm}\label{CPM}
Let $G_1\square G_2$ be the Cartesian product of two graphs $G_1$ and $G_2$. If $G_1\square G_2$ is connected, then the biharmonic distance between two vertices $u_1u_2$ and $v_1v_2$ is
$$(\sum_{1\leq i\leq n_1,1\leq j\leq n_2}(\lambda_{i}(G_1)+\lambda_{j}(G_2))^{-2}(z_i(G_1)\otimes z_j(G_2)(u_1u_2)-z_i(G_1)\otimes z_j(G_2)(v_1v_2))^2)^{\frac{1}{2}}.$$
\end{thm}
\subsection{The Cayley graph of a finite abelian group}
Let $\Gamma$ be a group and let $S$ be a subset of $\Gamma$ that is closed under taking inverses and does not contain the identity. The
{\em Cayley graph} $G(\Gamma, S)$, is the graph with vertex set $\Gamma$ and edge set $\{\{g,h\}\mid g^{-1}h\in S\}$. A
character $\chi$ of $\Gamma$ is a group homomorphism $\chi:\Gamma\rightarrow C^*$. Let $\Gamma^{\wedge}$ denote the set of characters of $\Gamma$. Denote the trivial character by $1_{\chi}$.
\begin{lemma}[\cite{Bab}]\label{L531}
Let $\Gamma$ be a finite abelian group. Then the adjacency eigenvalues of Cayley graph $G(\Gamma, S)$ are $\alpha_{\chi}$ with $\chi\in \Gamma^{\wedge}$, where $\alpha_{\chi}=\sum_{s\in S}\chi(s)$.
In addition, $\chi$ is an adjacency eigenvector with eigenvalue $\alpha_{\chi}$.
\end{lemma}

Note that the Cayley graph $G(\Gamma, S)$ is a regular graph and $L=|S|I-A$, where $L$ and $A$ are the Laplacian matrix and the adjacency matrix of $G(\Gamma, S)$ respectively. By Theorem \ref{m1} and Lemma \ref{L531}, we obtain the following result.
\begin{thm}\label{CAM}
Let $\Gamma$ be a finite abelian group. Then the biharmonic distance between two vertices $u$ and $v$ in Cayley graph $G(\Gamma, S)$ is
$$d_B(u,v)=(\sum_{1_{\chi}\neq\chi\in\Gamma^{\wedge}}(|S|-\alpha_{\chi})^{-2}(\chi(u)-\chi(v))^2)^{\frac{1}{2}}.$$
\end{thm}

\begin{figure}[hptb]
  \centering
  \includegraphics[width=8cm]{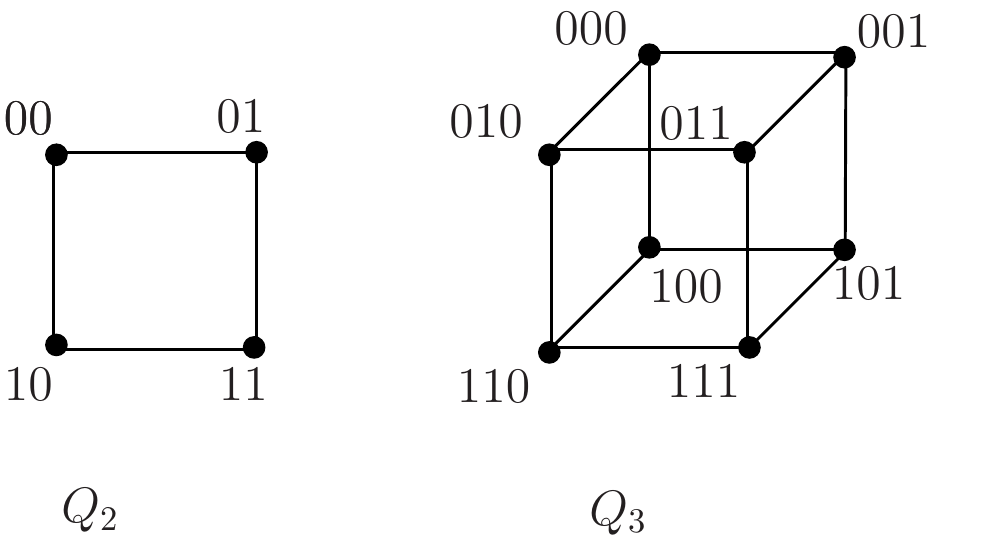}\\
  \caption{The structures of $Q_2$ and $Q_3$.
}\label{f2}
\end{figure}

The $n$-dimensional {\em hypercube} $Q_n$ \cite{Saa} is a graph with vertex set $\{x_{1}x_{2}\ldots x_{n}\mid x_i\in\{0, 1\}, 1\leq i\leq n\}$ and two vertices are adjacent if and only if they differ exactly in one position (see Figure \ref{f2}).

\begin{remark}\label{L532}
Let $(Z_2)^n$ be the direct products of the cyclic group $Z_2=\{0,1\}$ with $|(Z_2)^n|=2^n$. We use the notation $v_i$ to denote the $i$th coordinate of $v\in (Z_2)^n$, where $v=v_{1}v_{2}\ldots v_{n}$ and $v_k\in Z_2$ for $1\leq k\leq n$.
Note that the $n$-dimensional hypercube $Q_n$ can be viewed as a Cayley graph $G(\Gamma, S)$ satisfying that $\Gamma=(Z_2)^n$
and $S=\{e^i\mid 1\leq i\leq n\}$, where $e^i_h=1$ if and only if $h=i$. Since the characters of $(Z_2)^n$
are $\chi_I(v)=(-1)^{\sum_{i\in I} v_i}$ for $I\subseteq\{1,\ldots,n\}$, by Lemma \ref{L531}, the adjacency eigenvalue corresponding to $\chi_I$
is $\alpha_I=\sum_{j=1}^n\chi_I(e^j)=n-2|I|$.
\end{remark}

\begin{corollary}\label{cor1}
The biharmonic distance between two vertices $u$ and $v$ in $Q_n$ is
$$d_B(u,v)=\dfrac{\sqrt{2}}{2}(\sum_{\emptyset\neq I\subseteq\{1,\ldots,n\}}|I|^{-2}(1-(-1)^{\sum_{i\in I}(u_i+v_i)}))^{\frac{1}{2}}.$$
\end{corollary}
\proof By Theorem \ref{CAM} and Remark \ref{L532}, we have
\begin{eqnarray*}
d_B(u,v)&=&(\sum_{\emptyset\neq I\subseteq\{1,\ldots,n\}}(|S|-\alpha_I)^{-2}(\chi_I(u)-\chi_I(v))^2)^{\frac{1}{2}}\\
&=&(\sum_{\emptyset\neq I\subseteq\{1,\ldots,n\}}(n-(n-2|I|))^{-2}((-1)^{\sum_{i\in I}u_i}-(-1)^{\sum_{i\in I}v_i})^2)^{\frac{1}{2}}\\
&=&(\sum_{\emptyset\neq I\subseteq\{1,\ldots,n\}}(2|I|)^{-2}(2-2(-1)^{\sum_{i\in I}(u_i+v_i)}))^{\frac{1}{2}}\\
&=&\dfrac{\sqrt{2}}{2}(\sum_{\emptyset\neq I\subseteq\{1,\ldots,n\}}|I|^{-2}(1-(-1)^{\sum_{i\in I}(u_i+v_i)}))^{\frac{1}{2}}.
\end{eqnarray*}
This completes the proof of Corollary \ref{cor1}. $\qed$

\section{Conclusions}\label{F}
In this paper, we first give some characterizations of the biharmonic distance of a graph and obtain some bounds for it. Second, we reveal a relationship between biharmonic index and Kirchhoff index and determine the unique graph having the minimum biharmonic index among the connected graphs with given order. Finally, as applications, we investigate the biharmonic distance of some special graphs including the complement of a graph, the Cartesian product of two graphs and the Cayley graph of a finite abelian group. In particular, the combinatorial expression of biharmonic distance between any two vertices of hypercubes are determined.
Note that resistance distance has some interpretations based on random walks and electrical networks. A natural question is that whether the biharmonic distance has some similar interpretations. Thus, we state a few challenging open problems on the biharmonic distance of a graph.
\begin{pro}
Find a probability characterization of the biharmonic distance $d_B(u,v)$ between two vertices $u$ and $v$ for a graph $G$.
\end{pro}

\begin{pro}
Find a combinatorial (or physical) characterization of the biharmonic distance $d_B(u,v)$ between two vertices $u$ and $v$ for a graph $G$.
\end{pro}

\begin{pro}
Find lower and upper bounds on the biharmonic index $B(T)$ for trees $T$ of given order and characterize the extremal graphs.
\end{pro}
\section*{Acknowledgement}
R-H. Li's research is supported by the National Natural Science Foundation of China (No. 62072034 and 61772346). W. Yang's research is supported by the National Natural Science Foundation of China (No. 11671296).


\begin{thebibliography}{99}

\bibitem{Bab} L. Babai, Spectra of Cayley graphs, {\em J. Combin. Theory Ser. B}, 27(2) (1979), 180--189.  

\bibitem{Bam} B. Bamieh, M.R. Jovanovic, P. Mitra, S. Patterson, Coherence in large-scale networks:
Dimension-dependent limitations of local feedback, {\em IEEE Trans. Autom. Control}, 57(9) (2012), 2235--2249.  

\bibitem{Bap} R.B. Bapat, Graphs and matrices, Universitext, Springer, London, 2010.

\bibitem{Bar} P. Barooah, J.P. Hespanha, Graph effective resistances and distributed control: Spectral properties and applications, in {\em Proceedings of the 45th IEEE Conference on Decision and Control}, San Diego, CA, USA, 2006, 3479--3485.

\bibitem{Bon} J.A. Bondy, U.S.R. Murty, Graph Theory with Applications, The Macmillan Press Ltd, New York, 1976.

\bibitem{CheH} H. Chen, F. Zhang, Resistance distance and the normalized Laplacian spectrum,
{\em Discrete Appl. Math.}, 155(5) (2007), 654--661.  

\bibitem{CheD} X. Chen, K.C. Das, Characterization of extremal graphs from Laplacian eigenvalues and the sum of powers of the Laplacian eigenvalues of graphs, {\em Discrete Math.}, 338(7) (2015), 1252--1263.  

\bibitem{Cve} D. Cvetkovi\'{c}, P. Rowlinson, S. Simi\'{c}, 
An introduction to the theory of graph spectra, Cambridge University Press, Cambridge, 2010.

\bibitem{Das0} K.C. Das, A sharp upper bound for the number of spanning trees of a graph, {\em Graphs Combin.}, 23(6) (2007), 625--632.  

\bibitem{Das} K.C. Das, K. Xu, M. Liu, On sum of powers of the Laplacian eigenvalues of graphs, {\em Linear Algebra Appl.}, 439(11) (2013), 3561--3575. 

\bibitem{Doy} P.G. Doyle, J.L. Snell, Random walks and electric networks, Mathematical Association of America, Washington, DC, 1984.

\bibitem{Fie} M. Fiedler, Algebraic connectivity of graphs, {\em Czechoslovak Math. J.}, 23(2) (1973), 298--305.  

\bibitem{Fit} K. Fitch, N.E. Leonard, Joint centrality distinguishes optimal leaders in noisy networks, {\em IEEE Trans. Control Netw. Syst.}, 3(4) (2016), 366--378.  

\bibitem{Gut} I. Gutman, B. Mohar, The quasi-Wiener and the Kirchhoff indices coincide, {\em J. Chem. Inf. Comput. Sci.}, 36(5) (1996), 982--985. 

\bibitem{Kle} D.J. Klein, M. Randi\'{c}, Resistance distance, {\em J. Math. Chem.}, 12(1) (1993), 81--95.  

\bibitem{Lip} Y. Lipman, R.M. Rustamov, T.A. Funkhouser, Biharmonic distance, {\em ACM Transactions on Graphics}, 29(3) (2010), 1--11.  

\bibitem{Liu} M. Liu, B. Liu, A note on sum of powers of the Laplacian eigenvalues of graphs, {\em Appl. Math. Lett.}, 24(3) (2011), 249--252.  

\bibitem{Peng} Y.J. Peng, S.C. Li, On the Kirchhoff index and the number of spanning trees of linear phenylenes, {\em MATCH Commun. Math. Comput. Chem}, 77(3) (2017), 765--780.  

\bibitem{Saa} Y. Saad, M.H. Schultz, Topological properties of hypercubes, {\em IEEE Trans. Comput.}, 37(7) (1988), 867--872.  

\bibitem{Ste} W. Stevenson, Elements of Power System Analysis, McGraw Hill, New York, 1975. 

\bibitem{Xu} W. Xu, B. Wu, Z. Zhang, Z. Zhang, H. Kan, G. Chen, Coherence scaling of noisy second-order scale-free consensus networks, {\em IEEE Trans. Cybern.}, 52(7) (2022), 5923--5934.  

\bibitem{Yang} Y. Yang, D.J. Klein, A recursion formula for resistance distances and its applications, {\em Discrete Appl. Math.}, 161(16-17) (2013), 2702--2715.  

\bibitem{Yi} Y. Yi, B. Yang, Z. Zhang, S. Patterson, Biharmonic distance and the performance of second-order consensus networks with stochastic disturbances, in {\em 2018 Annual American Control Conference (ACC)}, Milwaukee, WI, USA, 2018, 4943--4950.

\bibitem{YiY} Y. Yi, B. Yang, Z. Zhang, Z. Zhang, S. Patterson, Biharmonic distance-based performance metric for second-order noisy consensus networks, {\em IEEE Trans. Inf. Theory}, 68(2) (2022), 1220--1236.  

\bibitem{Zhang} Z. Zhang, W. Xu, Y. Yi, Z. Zhang, Fast approximation of coherence for second-order noisy consensus networks, {\em IEEE Trans. Cybern.}, 52(1) (2022), 677--686.  

\bibitem{Zhou} B. Zhou, On sum of powers of the Laplacian eigenvalues of graphs, {\em Linear Algebra Appl.}, 429(8-9) (2008), 2239--2246.  
\end{thebibliography}
\end{document}